\documentclass[11pt]{amsart}
\usepackage{amssymb,amsmath,amsthm}
\usepackage{a4}
\usepackage{epsf}
\usepackage{epsfig}
\begin{document}

\overfullrule=3pt

\newcommand{\NP}{$\mathcal{N}\mathcal{P}$}
\newcommand{\newsymb}{{\mathcal P}}
\newcommand{\Pd}{{\mathcal P}^d}
\newcommand{\R}{{\mathbb R}}
\newcommand{\N}{{\mathbb N}}
\newcommand{\Q}{{\mathbb Q}}
\newcommand{\Z}{{\mathbb Z}}

\newcommand{\enorm}[1]{\Vert #1\Vert}
\newcommand{\inter}{\mathrm{int}}
\newcommand{\conv}{\mathrm{conv}}
\newcommand{\aff}{\mathrm{aff}}
\newcommand{\lin}{\mathrm{lin}}
\newcommand{\cone}{\mathrm{cone}}

\newcommand{\dist}{\mathrm{dist}}
\newcommand{\trans}{\intercal}
\newcommand{\diam}{\mathrm{diam}}
\newcommand{\vol}{\mathrm{vol}}

\newcommand{\pp}{\mathfrak{p}}
\newcommand{\pf}{\mathfrak{f}}
\newcommand{\pg}{\mathfrak{g}}
\newcommand{\PP}{\mathfrak{P}}
\newcommand{\pl}{\mathfrak{l}}
\newcommand{\pv}{\mathfrak{v}}
\newcommand{\cl}{\mathrm{cl}}
\newcommand{\bx}{\overline{x}}

\def\ip(#1,#2){#1\cdot#2}

\newtheorem{theorem}{Theorem}[section]
\newtheorem{corollary}[theorem]{Corollary}
\newtheorem{lemma}[theorem]{Lemma}
\newtheorem{remark}[theorem]{Remark}
\newtheorem{definition}[theorem]{Definition}
\newtheorem{proposition}[theorem]{Proposition}
\newtheorem{claim}[theorem]{Claim}
\numberwithin{equation}{section}

\title[Successive minima and Diophantine approximations]{Successive minima and best simultaneous Diophantine approximations }

\author{Iskander Aliev}
\address{Technische Universit\"at Wien,
Wiedner Hauptstraße 8-10 / 1046, 1040 Wien, \"Osterreich}
\email{ialiev@osiris.tuwien.ac.at}

\author{Martin Henk}
\address{Martin Henk, Universit\"at Magdeburg, Institut f\"ur Algebra und Geometrie,
  Universit\"ats\-platz 2, D-39106 Magdeburg, Germany}
\email{henk@math.uni-magdeburg.de}


\begin{abstract}
We study the problem of best approximations of a vector
$\alpha\in\R^n$ by rational vectors of a lattice $\Lambda\subset
\R^n$ whose common denominator is bounded. To this end we
introduce successive minima for a periodic lattice structure and
extend some classical results from geometry of numbers to this
structure. This leads to bounds for  the best approximation
problem which generalize and improve former results.

\end{abstract}

\keywords{successive minima, periodic lattice, simultaneous homogeneous approximation, critical determinant}
\subjclass[2000]{11J13, 11H06, 11H31, 52C07}

\maketitle

\section{Introduction}
Let $\mathcal{K}_0^n$ be the set of all $0$-symmetric convex bodies in
the $n$-dimensional Euclidean
space $\R^n$. For $K\in \mathcal{K}_0^n$ and $x\in\R^n$ we denote by $|x|_K = \min\{\lambda \geq 0 : x\in\lambda\,K\}$ the
norm of $x$ induced by $K$. If $K$ is the $n$-dimensional unit ball $B^n$ centered
at the origin then we write $\Vert x\Vert$ instead of $|x|_{B^n}$ and
the associated inner product is denoted by  $x\cdot y$, for  $x,y\in\R^n$.
 The volume, i.e., the $n$-dimensional Lebesgue
measure, of a set $S\subset\R^n$ is denoted by $\vol(S)$.

For a lattice $\Lambda=B\Z^n \subset \R^n$, $B\in \mathrm{GL(n,\R)}$,
a vector $\alpha\in\R^n$ and an integer $Q\geq 1$ the
functional
\begin{equation*}
     \beta(\alpha,Q,\Lambda,K):=\min\left\{ |q\,\alpha - b|_K >0 : q\in\{0,\dots,Q\},\, b\in\Lambda \right\}
\end{equation*}
measures the quality of a  best approximation of $\alpha$
by a {\em rational vector of the lattice
$\Lambda$} whose common denominator is bounded by $Q$.
This functional has been studied from various respects. For instance,
for $n=1$ and based on continued fractions, Klein \cite{klein:1895} gave a geometric
interpretation  of such a ``best
approximation point'' as a vertex of an associated $2$-dimensional
{\em Klein polyhedron}. Davenport and Mahler \cite{davenport-mahler:1946} proved that there exist
infinitely many points $(q,z)^\intercal\in\Z^3$, $z\in\Z^2$, such that
$\Vert q\alpha - z\Vert^2\leq (2/\sqrt{23})/q$ and the constant on
the right hand side is best possible. In the context of primal methods
for integer linear programs it has been also shown that
these best approximations are related to so called Hilbert Bases of
rational polyhedral pointed cones \cite{henk-weismantel:02}.

$\beta(\alpha,Q,\Lambda,K)$
was embedded by W.B.~Jurkat \cite{jurkat:81} and W. Kratz \cite{kratz:81}  in a series of functionals,
namely
for $1\leq i\leq n+1$ they defined
\begin{align*}
   \widetilde{\lambda}_i(\alpha,Q,\Lambda,K):=\min\big\{ \lambda\geq 0 : &\,\,
   \exists \text{  $i$ linearly independent points }(q_j, b_j)^\intercal\in\R^{n+1}, \\
    &\,q_j\in\{0,\dots,Q\}, b_j\in\Lambda\text{ with } \vert q_j\alpha -
    b_j\vert_K \leq \lambda \big\}.
\end{align*}
$\widetilde{\lambda}_i(\alpha,Q,\Lambda,K)$ is called the
$i$-th successive minimum with respect to $\alpha, Q, \Lambda$ and $K$.
For abbreviation we just write $\widetilde{\lambda}_i$, since the
dependency on $\alpha$ etc.~will be clear from the
context.
Note that $\widetilde{\lambda}_1=\beta(\alpha,Q,\Lambda,K)$ and
$\widetilde{\lambda}_i\leq \widetilde{\lambda}_{i+1}$.

In \cite{kratz:81} successive-minima-type inequalities with
constraints were studied and, among others, the following inequalities were
proven in the case $K=B^n$
\begin{equation}
\label{eq:lower_and_upper}
 \frac{2^{\frac{n+1}{2}}}{(n+1)!}\det\Lambda \leq
 \widetilde{\lambda}_1\cdot\widetilde{\lambda}_2\cdots \widetilde{\lambda}_{n+1}\,\vol(B^{n+1})\,\max_{1\leq i\leq n+1}\left\{\frac{q_i}{\widetilde{\lambda}_i}\right\}\leq 2^{2n+1}\det\Lambda,
\end{equation}
where
 $q_i$  belongs to a point $(q_i,b_i)^\intercal$ attaining the
$i$-th successive minimum, i.e., $\widetilde{\lambda}_i=\Vert
q_i\alpha-b_i\Vert$.
 For $n=2$
these inequalities  were improved by Kratz \cite{kratz:99} to
\begin{equation}
\label{eq:twodimall}
  \frac{2}{3\,\sqrt{3}} \det\Lambda \leq\widetilde{\lambda}_1\,\widetilde{\lambda}_2\,\widetilde{\lambda}_{3}\,\max_{1\leq i\leq 3}\left\{\frac{q_i}{\widetilde{\lambda}_i}\right\} \leq \frac{2}{\sqrt{3}}\det\Lambda
\end{equation}
and moreover it was shown
\begin{equation}
\label{eq:twodim}
\widetilde{\lambda}_1\,\widetilde{\lambda}_2 \leq\frac{2}{\sqrt{3}}\frac{\det\Lambda}{Q}.
\end{equation}
All the constants in \eqref{eq:twodimall} and \eqref{eq:twodim} are
best possible.
Inequalities of that type give us information on the quality of the
simultaneous approximation of a vector by  a system of  rational
vectors of a lattice whose common denominators are bounded.

The inequalities \eqref{eq:lower_and_upper} and  \eqref{eq:twodimall} may be regarded
as analogs (in the case $K=B^n$) to
Minkowski's classical inequalities on successive minima
(cf.~\cite[pp.~59]{gruber-lekkerkerker:87})
\begin{equation}
\label{eq:minkowski_second}
  \frac{2^n}{n!} \det\Lambda \leq \lambda_1(\Lambda,K)\cdot\ldots\cdot\lambda_n(\Lambda,K) \vol(K) \leq 2^n\,\det\Lambda.
\end{equation}
Here the $i$-th successive
minimum $\lambda_i(\Lambda,K)$ is defined as
\begin{equation*}
    \lambda_i(\Lambda,K):=\min\{\lambda \geq 0: \dim(\Lambda\cap\lambda\,K)\geq i\}.
\end{equation*}
Both bounds in \eqref{eq:minkowski_second} are best possible.
Statement \eqref{eq:twodim} is the $2$-dimensional analog of another result of
Minkowski (cf.~\cite[pp.~195]{gruber-lekkerkerker:87}) on successive minima of $n$-dimensional unit ball
\begin{equation}
\label{eq:minkowski_ball}
  \lambda_1(\Lambda,B^n)\cdot\ldots\cdot\lambda_n(\Lambda,B^n) \Delta(B^n) \leq \det\Lambda.
\end{equation}
Here $\Delta(K)$ denotes the critical determinant of $K\in\mathcal{K}_0^n$,
i.e., the minimal determinant of a lattice whose only lattice point
belonging to the interior of $K$ is the origin. Since $\Delta(B^2)=\frac{\sqrt{3}}{2}$ (cf.~\cite[pp.~244]{gruber-lekkerkerker:87}, \cite[pp.~8]{zong:99a}) the analogy between \eqref{eq:twodim} and \eqref{eq:minkowski_ball} is obvious.

The fact that the constant $\frac{2}{\sqrt{3}}$ in \eqref{eq:twodim} is best
possible  follows also from a recent and more general
result by I.~Aliev and P.~M.~Gruber \cite{aliev-gruber:04}, who showed that for every
$\epsilon>0$ there exists a vector $\alpha\in\R^n\setminus\{0\}$ and a $Q\in\N$
such that
\begin{equation}
\label{eq:aliev_gruber}
       \left(\widetilde{\lambda}_1 \right)^n>
       \frac{1-\epsilon}{\Delta(K)} \frac{\det\Lambda}{Q}.
\end{equation}
Hence \eqref{eq:twodim} cannot be improved. In fact,
\eqref{eq:aliev_gruber} was not only proven for $0$-symmetric
convex bodies, but for any bounded star body $K$. Finally, it is
also known (cf.~\cite[p.~197]{gruber-lekkerkerker:87}) that in the
planar case \eqref{eq:minkowski_ball} can be generalized to all
$K\in\mathcal{K}_0^n$, i.e.,
\begin{equation}
 \lambda_1(\Lambda,K)\cdot\lambda_2(\Lambda,K) \Delta(K) \leq \det\Lambda.
\label{eq:minkowski_planar}
\end{equation}

In this paper  we want  to introduce and study a slightly different series
of successive minima associated to this best approximation
problem. With these successive minima we can give best possible upper
and lower bounds of a similar type as in \eqref{eq:lower_and_upper}
and \eqref{eq:twodim} with respect to all
$0$-symmetric convex bodies and, in addition, these results include the
classical inequalities \eqref{eq:minkowski_second},
\eqref{eq:minkowski_ball} and \eqref{eq:minkowski_planar}.
To this end we consider the special periodic lattice
\begin{equation*}
     \Lambda(\alpha,Q):=\Lambda\cup(\alpha+\Lambda)\cup(2\,\alpha+\Lambda)\cup\dots\cup  (Q\,\alpha+\Lambda),
\end{equation*}
where we always assume that $k\alpha\notin\Lambda$ for $1\leq k\leq
Q$. In order to allow the case $\alpha\in\Lambda$ we admit $Q=0$.
Next we define for $1\leq i\leq n$
\begin{equation*}
   \lambda_i(\Lambda(\alpha,Q),K):=\min\{\lambda\geq 0 :
   \dim(\Lambda(\alpha,Q)\cap\lambda\,K)\geq i\},
\end{equation*}
and again for abbreviation we just write $\lambda_i$ instead of
$\lambda_i(\Lambda(\alpha,Q),K)$. In comparison with the
successive minima $\widetilde{\lambda}_i$ we note that
\begin{equation*}
  \widetilde{\lambda}_1 = \lambda_1,\quad\text{ and } \quad
  \widetilde{\lambda}_i\leq \lambda_i,\quad 2\leq i\leq n.
\end{equation*}
Hence any upper bound on $\lambda_i$ gives us also an upper bound on
$\widetilde{\lambda}_i$.
In order to state our results concerning these successive minima we
need some more notation. For $K\in\mathcal{K}_0^n$ we denote by
 $\delta(K)$ the density of a
densest packing of translates of $K$
(cf.~\cite[pp.~218]{gruber-lekkerkerker:87}) and the dual lattice of
$\Lambda$ is denoted by $\Lambda^*$ (cf.~\cite[pp.~23]{gruber-lekkerkerker:87}).

\begin{theorem} Let $K\in\mathcal{K}_0^n$, $\alpha\in\R^n$,
                $\Lambda\subset\R^n$ be a lattice and $Q\in\N_{\geq 0}$
                such that $k\alpha\notin \Lambda$ for $1\leq k \leq
                Q$. Then with
                $\lambda_i=\lambda_i(\Lambda(\alpha,Q),K)$ we have
\begin{align*}
 \text{\rm i)}&\quad  (\lambda_1)^n \vol(K) \leq
 \delta(K)\,2^n\,\frac{\det\Lambda}{Q+1}, \\
\text{\rm ii)} &\quad \frac{2^n}{n!}
\det\Lambda\,\gamma(\alpha,\Lambda,Q,n)\leq
\lambda_1\cdot\ldots\cdot\lambda_n\, \vol(K) \leq
2^n\,\frac{\det\Lambda}{Q+1},
\end{align*}
where
\begin{align*}
  \gamma(\alpha,\Lambda,Q,n)=\min\Big\{ &|u^*\cdot\alpha+z|>0 :  u^*\in \Lambda^*,\,z\in\Z, \\
 &\,\,\Vert (u^*,z)^\intercal\Vert \leq (n\lambda_n+Q\Vert(\alpha,1)^\intercal\Vert)^n/\det\Lambda\Big\}.
\end{align*}
\label{thm:minkowski}
\end{theorem}
\begin{remark} \hfill
\begin{enumerate} {\rm
\item If $K=[-1,1]^n$ is the standard cube of edge length 2
  centered at the origin, i.e., $|\cdot|_K$ is the maximum norm,
  and if $\Lambda=\Z^n$ statement i) implies that for
  any $(\alpha_1,\dots,\alpha_n)^\intercal\in\R^n$ there exist a
  $(z_1,\dots,z_n)^\intercal\in\Z^n$ and $q\in\{1,\dots,Q\}$ such that
\begin{equation*}
        |q\alpha_i-z_i| < Q^{-1/n},\quad 1\leq i\leq n.
\end{equation*}
 This is  Dirichlet's classical approximation theorem.
\item If $(Q+1)\,\alpha\in\Lambda$ then $\Lambda(\alpha,Q),K)$ is a
  lattice of determinant $\det\Lambda/(Q+1)$  and we also have
  $\gamma(\alpha,\Lambda,Q,n) \geq 1/(Q+1)$. Thus, in this situation,
  statement ii) becomes \eqref{eq:minkowski_second}. In particular,
  these inequalities are best possible.
\item In general we cannot expect to find a lower bound in ii) which
  depends only on $\Lambda$, $n$ and $Q$. For instance, let
  $\Lambda=\Z^n$ and for a positive integer $m$ let $K(m)$ be the
  cross-polytope with vertices $\{\pm \frac{1}{m}e_1,\pm e_i: 2\leq
  i\leq n\}$.  Here $e_i\in\R^n$ denotes the $i$-th unit vector. Then  $\vol(K(m))=(1/m)2^n/n!$ and for
  $\alpha=\frac{1}{m}e_1$, $Q<m$, we  have $\lambda_i=1$, $1\leq
  i\leq n$, and so  $\lambda_1\cdot\ldots\cdot\lambda_n\,\vol(K)=
  (1/m) 2^n/n!$. In this case we also have
  $\gamma(\alpha,\Z^n,Q,n) = 1/m$ and thus equality in the lower bound.
}
\end{enumerate}
\end{remark}

For the special cases $K=B^n$ or $n=2$ we obtain the following improvements on the upper bound in Theorem \ref{thm:minkowski} ii).
\begin{proposition} With the notation as  in Theorem \ref{thm:minkowski}
  we have
\begin{align*}
\text{\rm i)}&\quad
\lambda_1\cdot\lambda_2\cdot\ldots\cdot\lambda_n\vol(B^n)\leq
\delta(B^n)\, 2^n \, \frac{\det\Lambda}{Q+1}, \\
\text{\rm ii)}&\quad \lambda_1\cdot\lambda_2\,\Delta(K)\leq
\frac{\det\Lambda}{Q+1},\quad\text{ for }n=2.
\end{align*}
\label{prop:minkowski}
\end{proposition}
\begin{remark} \hfill\begin{enumerate}
 \item {\rm By \eqref{eq:aliev_gruber} we see that ii) is best
     possible and it generalizes \eqref{eq:twodim} to all norms in
     $\R^2$.
  \item Since $\Delta(K)= \vol(K)\,2^{-n}/\delta_L(K)$ where
   $\delta_L(K)$ denotes the density of a densest lattice packing of $K$
   (cf.~\cite[p.~221]{gruber-lekkerkerker:87}),
   inequality i) is a bit weaker than
   \eqref{eq:minkowski_ball}. However, since we deal with more general
   structures than lattice we have to replace, in comparison with
   \eqref{eq:minkowski_ball}, $\delta_L(K)$ by $\delta(K)$.
\item It is still an open conjecture of Davenport that
  \eqref{eq:minkowski_ball} can be generalized to arbitrary
  $K\in\mathcal{K}_0^n$ and  the same can be conjectured for inequality i).
}
\end{enumerate}
\end{remark}
According to the announced proofs of the
Kepler-Conjecture (cf.~\cite{hales:web}, \cite{hsiang:01}) we also
know $\delta(B^3)=\delta_L(B^3)=\pi/\sqrt{18}$ and therefore,  Proposition
\ref{prop:minkowski} i) leads to
\begin{corollary} For $K=B^3$ we have
\begin{equation*}
\label{eq:3dimball}
   \lambda_1\cdot\lambda_2\cdot \lambda_3\leq
\sqrt{2}\,\frac{\det\Lambda}{Q+1}.
\end{equation*}
\end{corollary}
Again, from \eqref{eq:aliev_gruber} we see that this inequality
is best possible.


\section{Proofs}
The basic properties of our set $\Lambda(\alpha,Q)$,
which allow us to extend inequalities on
lattices and convex bodies to this structure  are
\begin{equation}
\begin{split}
 {\rm i)} &\quad \Lambda(\alpha,Q) \text{ is invariant with respect to
   lattice translations of } \Lambda. \\
 {\rm ii)} &\quad \text{A fundamental cell of } \Lambda  \text{
   contains exactly }
   Q+1 \text{ points of } \Lambda(\alpha,Q).\\
 {\rm iii)} & \quad \text{For } a_1,a_2\in \Lambda(\alpha,Q) \text{ at least one of the points }  a_1-a_2,\; a_2-a_1\\
    &\quad\text{belongs to }\Lambda(\alpha,Q).
\end{split}
\label{eq:properties}
\end{equation}
A  set $S$ satisfying the first two properties is called a periodic
lattice. For such a periodic lattice and a set $X$ the arrangement $S+X$
is called a periodic lattice packing if for all $c_1, c_2\in S$, $c_1\ne
c_2$,
\begin{equation*}
     c_1+\inter(X)\cap c_2+\inter(X) =\emptyset,
\end{equation*}
where $\inter(\cdot)$ denotes the interior.
If $X$ is measurable and bounded then the density $\delta(S,X)$
of such a periodic packing can be calculated by (cf.~\cite[pp.~26]{rogers:64})
\begin{equation}
   \delta(S,X) = \vol(X)\frac{Q+1}{\det\Lambda}.
\label{eq:density}
\end{equation}
Now let $\lambda_1=\lambda_1(\Lambda(\alpha,Q),K)$. First we note
that on account of \eqref{eq:properties} iii)
\begin{equation*}
   \lambda_1=\min\left\{|a_1-a_2|_K : a_1,a_2\in
   \Lambda(\alpha,Q)), a_1\ne a_2\right\}.
\label{eq:scale_property}
\end{equation*}
Hence  $\frac{2}{\lambda_1}\Lambda(\alpha,Q)+K$ is a periodic lattice packing of $K$ and thus
\begin{equation*}
  \delta(K)\geq \delta\left(\frac{2}{\lambda_1}\Lambda(\alpha,Q),K\right) =
  (Q+1)\frac{\vol(K)}{\det\left(\frac{2}{\lambda_1}\Lambda\right)} =
  (Q+1) \frac{\lambda_1^n}{2^n} \frac{\vol(K)}{\det\Lambda}.
\label{eq:proofi}
\end{equation*}
This shows already Theorem \ref{thm:minkowski} i).

As an immediate consequence we have the following analog to a
theorem of Blichfeldt  \cite[p.~42]{gruber-lekkerkerker:87}
\begin{lemma} Let $X\subset \R^n$ be a measurable  set
  with $\vol(X) >  \det(\Lambda)/(Q+1)$. Then there exist two distinct
  $x_1,x_2\in X$ such that $x_1-x_2 \in \Lambda(\alpha,Q)$.
\label{lem:blichfeldt}
\end{lemma}
\begin{proof} W.l.o.g.~let $X$ be bounded. Now
  suppose the contrary, i.e., for all $x_1,x_2\in
  X$, $x_1\ne x_2$,  we have  $x_1-x_2 \notin \Lambda(\alpha,Q)$. Then $\Lambda(\alpha,Q)+X$
  is a periodic lattice packing of $X$, because $x\in
  c_i+\inter(X)$ for two distinct $c_1,c_2\in
  \Lambda(\alpha,Q)$ implies $x-c_1,x-c_2\in X$. Together with  \eqref{eq:properties} iii) this shows that $x-c_1-(x-c_2)=c_2-c_1\in\Lambda(\alpha,Q)$
  or $x-c_2-(x-c_1)=c_1-c_2\in\Lambda(\alpha,Q)$.

  Thus $\Lambda(\alpha,Q)+X$
  is a periodic lattice packing of $X$ and  by
  \eqref{eq:density} we obtain the contradiction
\begin{equation*}
     1\geq \delta(\Lambda(\alpha,Q),X)= \vol(X)\,\frac{Q+1}{\det\Lambda}.
\end{equation*}
\end{proof}

Now we come to
\begin{proof}[Proof of Theorem \eqref{thm:minkowski} {\rm ii)}] On
  account of Lemma \ref{lem:blichfeldt} the
  upper bound can be shown completely analogous to the upper bound  of
  Minkowski's classical inequality \eqref{eq:minkowski_second} . 
   Here we will follow
  Siegel's proof \cite[pp.33]{siegel:89} and we will just give the main
  arguments.

  Let $a_1,\dots,a_n\in \Lambda(\alpha,Q)$ be linearly independent such that $a_i\in
  \lambda_i\,K$, $1\leq i\leq n$.  For $x\in\R^n$ and $0\leq k\leq
  n-1$ we denote by $L_k(x)$ the $k$-dimensional affine plane given by
  $L_k(x)=x+\lin\{a_1,\dots,a_k\}$, where $\lin\{\}$ denotes the
  linear hull. Let $\widetilde{K}=\inter(K)$ and for $x\in \widetilde{K}$ let $c_k(x)$ be the center of gravity of
  the intersection $L_k(x)\cap \widetilde{K}$. Since $\widetilde{K}$ is convex $c_k(x)$
  belongs to $\widetilde{K}$, and moreover, $c_k(x)$ depends continuously on $x$.
  Let $f:\widetilde{K}\to\R^n$ be the map
\begin{equation}
    f(x)=\lambda_1\,c_0(x)+(\lambda_2-\lambda_1)\,c_1(x)+\cdots + (\lambda_{n}-\lambda_{n-1})\,c_{n-1}(x).
\label{eq:map}
\end{equation}
Then it is shown that $f$ is an injective map and that the volume of the
bounded set $C=f(\widetilde{K})$ (which, in general,  is not convex) is given by (cf.~\cite[pp.33]{siegel:89})
\begin{equation}
   \vol(C)=\lambda_1\cdot\lambda_2\cdot\ldots\cdot\lambda_n\vol(K).
\label{eq:vol}
\end{equation}
Next we claim
\begin{equation}
  \Lambda(\alpha,Q) \text{ is a periodic lattice packing of }\frac{1}{2}C.
\label{eq:claim}
\end{equation}
Suppose the opposite and let
$y_1,y_2\in C$, $y_1\ne y_2$,  such that
$\frac{1}{2}(y_1-y_2)\in\Lambda(\alpha,Q)$.
Let $x_1,x_2\in \widetilde{K}$ with $y_i=f(x_i)$ and let $r$ be minimal such
that $x_1-x_2\in\lin\{a_1,\dots,a_r\}$. Then $c_k(x_1)=c_k(x_2)$ for
$k=r,\dots,n-1$ and therefore we may write
\begin{equation}
\begin{split}
 \frac{1}{2}(y_1-y_2)=&\lambda_1\frac{1}{2}\left(c_0(x_1)-c_0(x_2)\right)+
(\lambda_2-\lambda_1)\frac{1}{2}\left(c_1(x_1)-c_1(x_2)\right) +\\
&\cdots +
(\lambda_{r}-\lambda_{r-1})\frac{1}{2}\left(c_{r-1}(x_1)-c_{r-1}(x_2)\right).
\end{split}
\label{eq:strange}
\end{equation}
Since $\widetilde{K}$ is a convex $0$-symmetric open set we have
$\frac{1}{2}\left(c_i(x_1)-c_i(x_2)\right)\in \widetilde{K}$ and so we
find by \eqref{eq:strange}
\begin{equation}
\frac{1}{2}(y_1-y_2)\in \inter\{\lambda_r\,K\}.
\label{eq:contradict}
\end{equation}
On the other hand, \eqref{eq:strange} also shows
$\frac{1}{2}(y_1-y_2)\in
\lambda_r\frac{1}{2}(x_1-x_2)+\lin\{a_1,\dots,a_{r-1}\}$ and thus, by
the choice of $r$, the point
$\frac{1}{2}(y_1-y_2)$ is linearly independent from
$a_1,\dots,a_{r-1}$.  In view of \eqref{eq:contradict}, however, this
contradicts the definition of $\lambda_r$.

Thus \eqref{eq:claim} is shown and by Lemma
\ref{lem:blichfeldt} (applied to $\frac{1}{2}C$) and \eqref{eq:vol} we obtain
the upper bound in Theorem \ref{thm:minkowski} ii).

For the lower bound let $a_i$ as above. Then $K$ contains the cross-polytope with vertices $\{\pm \frac{1}{\lambda_i}a_i : 1\leq i\leq n\}$. Hence
\begin{equation}
   \vol(K)\geq
   \frac{2^n}{n!}\frac{1}{\lambda_1\cdot\ldots\cdot\lambda_n}\vert\det(a_1\,a_2\,\dots\,a_n)\vert.
\label{eq:in1}
\end{equation}
With $a_i=b_i-q_i\,\alpha$, $b_i\in\Lambda$, $q_i\in\{0,\dots,Q\}$, $1\leq i\leq n$,  we may write
\begin{equation}
  \det(a_1\,a_2\,\dots\,a_n) =\det\begin{pmatrix} b_1 & b_2 & \cdots & b_n & \alpha \\
    q_1 & q_2 & \cdots & q_n & 1
\end{pmatrix} =:\det A.
\label{eq:in2}
\end{equation}
Now let ${\Lambda^\prime}=\{(b,z)^\intercal\in\R^{n+1} :
b\in\Lambda, z\in\Z \}=\Lambda\times \Z$ be the lattice given by the Cartesian product of
$\Lambda$ and $\Z$.  Obviously, we have ${\Lambda^\prime}^\star =
\Lambda^\star \times \Z$. Let ${u^\prime}^*=(u^*,z)\in
{\Lambda^\prime}^\star$ be the uniquely (up to $\pm$) determined
primitive vector which is orthogonal  to the lattice hyperplane $L$
of $\Lambda^\prime$
generated by $(b_1,q_1)^\intercal,\dots, (b_1,q_n)^\intercal$. Next we
supplement ${u^\prime}^*$ to a basis of ${\Lambda^\prime}^\star$ and
let $A^*$ be the matrix consisting of these basis vectors as column
vectors. Since the inner product of a vector of  lattice with a vector
of the dual lattice is an integer we get
\begin{equation}
\det(A) = \frac{1}{\det{\Lambda^\prime}^\star}\det(A^\intercal A^*)
\geq  \det\Lambda\,   \vert u^*\alpha+z\vert.
\label{eq:in3}
\end{equation}
On the other hand we have (cf. \cite[pp.~28]{martinet:03})
\begin{equation*}
   \det(L\cap \Lambda^\prime)=\det\Lambda^\prime\cdot\Vert {u^\prime}^*\Vert=\det\Lambda\,\Vert{u^\prime}^*\Vert,
\end{equation*}
and since the $n$ linearly independent points $(b_i,q_i)^\intercal$ belong
to $L\cap\Lambda^\prime$ we can bound $\Vert{u^\prime}^*\Vert$ by
\begin{equation}
   \Vert{u^\prime}^*\Vert \leq \frac{\det(L\cap{\Lambda^\prime})}{\det\Lambda}\leq
   \frac{\Vert(b_1,q_1)^\intercal\Vert\cdot\ldots\cdot
     \Vert(b_n,q_n)^\intercal\Vert}{\det\Lambda}.
\label{eq:in4}
\end{equation}
Finally we observe that
\begin{equation*}
 \Vert(b_i,q_i)^\intercal \Vert\leq \Vert (b_i,q_i)^\intercal - q_i
 (\alpha,1)^\intercal\Vert + Q\Vert   (\alpha,1)^\intercal\Vert \leq \sqrt{n}\,\lambda_i+Q\Vert   (\alpha,1)^\intercal\Vert
\end{equation*}
and together with \eqref{eq:in1}--\eqref{eq:in4} we get the desired
lower bound.
\end{proof}

\begin{proof}[Proof of Proposition \ref{prop:minkowski}]
As in the case of \eqref{eq:minkowski_ball} and  \eqref{eq:minkowski_planar} with respect to Minkowski's successive minima, these
inequalities are to some extend just a consequence of the proof
of the upper bound in Theorem \ref{thm:minkowski} ii). Roughly
speaking, the main observation is that for $K=B^n$ or $n=2$ the
map $f$ in \eqref{eq:map} can be made linear. For details we
refer to \cite[pp.~195, p.~197]{gruber-lekkerkerker:87}.

If $f$ is linear the resulting body $f(K)$ is still a $0$-symmetric
convex body and from \eqref{eq:claim} we know that $\Lambda(\alpha,Q)$ is a periodic lattice packing  of
  $\frac{1}{2}f(K)$. Thus we know that $\lambda_1(\Lambda(\alpha,Q),f(K))\geq
  1$ and with Theorem   \ref{thm:minkowski} i) we get for $K=B^n$ or
  $n=2$
\begin{equation*}
   \lambda_1\cdot\ldots\cdot\lambda_n\,\vol(K)=\vol(f(K)) \leq
   \delta(f(K)) 2^n \frac{\det\Lambda}{Q+1} = \delta(K) 2^n \frac{\det\Lambda}{Q+1},
\end{equation*}
since the density of a densest packing is invariant with respect to
affine transformations.  For $K=B^n$ we get Proposition
\ref{prop:minkowski} i). In the planar case it is well known that
$\delta(K)=\delta_L(K)= \vol(K)\,2^{-n}/\Delta(K)$
(cf.~\cite[pp.~248]{gruber-lekkerkerker:87}) and so we also obtain the second statement.
\end{proof}

\section{Acknowledgement}

The first author was supported by FWF Austrian Science Fund, projects M672 and M821-N12.

\bibliographystyle{plain}
\def\cprime{$'$} \def\cprime{$'$} \def\cprime{$'$} \def\cprime{$'$}
  \def\cprime{$'$} \def\cprime{$'$} \def\cprime{$'$}
\providecommand{\bysame}{\leavevmode\hbox to3em{\hrulefill}\thinspace}
\providecommand{\MR}{\relax\ifhmode\unskip\space\fi MR }
\providecommand{\MRhref}[2]{%
  \href{http://www.ams.org/mathscinet-getitem?mr=#1}{#2}
}
\providecommand{\href}[2]{#2}

\end{document}